\documentclass[11pt]{article}
\usepackage{amsmath,amsthm,amsfonts,amssymb,amscd}
\usepackage{bbm}
\usepackage{enumerate}
\usepackage{enumitem}
\usepackage{url}
\usepackage{cite}
\usepackage{fullpage}
\usepackage{fancyhdr}
\usepackage{amsmath,amsthm,amsfonts,amssymb,amscd}
\usepackage{color}
\usepackage{float}
\usepackage{soul}
\usepackage{graphicx}
\usepackage[margin=1.0in]{geometry}
\theoremstyle{plain}
\newtheorem{theorem}{Theorem}[section]
\newtheorem{lemma}{Lemma}[section]
\newtheorem{definition}{Definition}[section]

\newtheorem{proposition}{Proposition}[section]
\newtheorem{example}{Example}[section]

\theoremstyle{definition}

\newcommand{\la}{\langle}
\newcommand{\ra}{\rangle}

\newcommand{\nexto}{\kern -0.54em}

\newcommand{\dZ}{{\cal Z \kern -0.7em Z}}
\newcommand{\dC}{{\rm\hbox{C \kern-0.8em\raise0.2ex\hbox{\vrule height5.4pt width0.7pt}}}}
\newcommand{\dQ}{{\rm\hbox{Q \kern-0.85em\raise0.25ex\hbox{\vrule height5.4pt width0.7pt}}}}

\usepackage{epsfig}
\usepackage{latexsym}

\newcommand{\NN}{\mathbb{N}}
\newcommand{\HH}{\mathcal{H}}
\newcommand{\RR}{\mathbbm{R}}

\newcommand{\dsty}{\displaystyle}

\newenvironment{retraitsimple}{\begin{list}{--~}{
 \topsep=0.3ex \itemsep=0.3ex \labelsep=0em \parsep=0em
 \listparindent=1em \itemindent=0em
 \settowidth{\labelwidth}{--~} \leftmargin=\labelwidth
}}{\end{list}}

\begin{document}

\title{On weak and strong convergence of the projected gradient method for convex optimization in real Hilbert spaces}
\author{J.Y. Bello Cruz\footnote{Instituto de Matem\'atica e Estat\'istica, Universidade Federal de
Goi\'as. Campus Samambaia, CEP 74001-970, Goi\^ania - G.O., Brazil.
E-mail: yunier@impa.br $\&$ yunier.bello@gmail.com}\and W. de
Oliveira\footnote{Universidade do Estado do Rio de Janeiro. E-mail: welington@ime.uerj.br}} 
\maketitle
\begin{abstract}
\noindent This work focuses on convergence analysis of the projected
gradient method for solving constrained convex minimization problems in Hilbert spaces. We show
that the sequence of points generated by the method employing the
Armijo line search converges weakly to a solution of the considered
convex optimization problem. Weak convergence is established by
assuming convexity and Gateaux differentiability of the objective
function, whose Gateaux derivative is supposed to be uniformly
continuous on  bounded sets. Furthermore, we propose some
modifications in the classical projected  gradient method in order
to obtain strong convergence. The new variant has the following
desirable properties: the sequence of generated points  is entirely
contained in a ball with diameter equal to the distance between the
initial point and the solution set, and the whole sequence converges
strongly to the solution of the problem that lies closest to the
initial iterate. Convergence analysis of both methods is presented without Lipschitz continuity assumption.
\end{abstract}

{\small
\noindent{\bf Keywords:} Armijo line search; Convex minimization;
Projection method, Strong and weak convergence.


\noindent{\bf Mathematical Subject Classification (2010):} 90C25,
90C30.
}
\section{Introduction}
In this work we are interested in weak and strong convergence of
projected gradient methods for the following optimization problem
\begin{equation}\label{1}
\min \,f(x) \;\mbox{ s.t. } \;x\in C\subset \HH\,,
\end{equation} under the following assumptions:
$C$ is a nonempty, closed and convex subset of a real Hilbert space $\HH$;
$f:\HH\rightarrow \RR$ is a convex and Gateaux differentiable function, and  $S_*$ is the solution set of problem~\eqref{1}, which is assumed to be nonempty.
Moreover, we assume the Gateaux derivative $\nabla f$ of $f$ to be uniformly continuous on bounded sets. The latter is
a weaker assumption than the commonly adopted one, which requires $\nabla f$ to be Lipschitz continuous in the whole space $\HH$.

Due to its simplicity, the projected gradient method has been widely
used in practical applications. The method has several useful
advantages. Primarily, it is easy to implement (especially, for
optimization problems with relatively simple constraints). The
method uses little storage and readily exploits any sparsity or
separable structure of $\nabla f$ or $C$. Furthermore, it is able to
drop or add active constraints during the iterations. Some important
references on the projected gradient method in finite dimensional
spaces are, for instance, \cite{polyak_libro, bertsekas}.

A general description of the classical projected gradient method can be stated as follows,
where we denote the projection of a given point $x$ onto $C$ by $P_C(x)$.
\begin{center}\fbox{\begin{minipage}[b]{\textwidth}
\noindent {\bf Projected Gradient Method}

\medskip

\noindent{\bf Initialization Step:} Take $x^0\in C$ and set $k=0$.

\noindent{\bf Iterative Step:} Given $x^k$, compute
\begin{equation}\label{9}
z^k = x^k -\beta_k \nabla f(x^k)
\end{equation}
\begin{equation}\label{10}
x^{k+1} = \alpha_k P_C(z^k) + (1-\alpha_k) x^k,
\end{equation}
where $\alpha_k\in (0,1]$ and $\beta_k$ is positive for all $k$. 

\noindent{\bf Stop Criterion:} If $x^{k+1}=x^k$ then stop.

\noindent Otherwise, set $k=k+1$ and go back to {\bf Iterative Step}.
\end{minipage}}\end{center}
Several choices are possible for
the stepsizes $\beta_k$ and $\alpha_k$. We focus our attention on the description of the following four strategies:

\begin{retraitsimple}

\item [{\bf(a)} ] \emph{Constant stepsize}: $\beta_k=\beta$ for all $k$ where $\beta>0$ is a
fixed number and $\alpha_k=1$ for all $k$.

\item [{\bf(b)} ] \emph{Armijo line search along the boundary of $C$}:  $\alpha_k = 1$
for all $k$ and $\beta_k$ is determined by
$
\beta_k=\bar{\beta}\theta^{\ell(k)}
$, 
for some $\bar{\beta}>0$, $\theta, \delta\in (0,1)$ where
\[
\ell(k):=\min\left\{\ell\in\NN \mid  f(P_C(x^{k,\ell}))\le
f(x^k)-\delta\la\nabla f(x^k), x^k-P_C(x^{k,\ell})\ra\right\}, \;
x^{k,\ell}=x^k-\bar{\beta}\theta^{\ell}\nabla f(x^k).
\]
\item [{\bf(c)} ] \emph{Armijo line search along the feasible direction}:
$(\beta_k)_{k\in \NN}\subset[\check{\beta},\hat{\beta}]$ for some
$0<\check{\beta}\le\hat{\beta}<\infty$ and $\alpha_k$ determined by
the following Armijo rule $ \alpha_k=\theta^{j(k)}$,  for some
$\theta, \delta\in (0,1)$ where
\[
j(k):=\min \left\{j\in \NN \mid f(x^{k,j})\le f(x^k)-\delta
\theta^{j}\la\nabla f(x^k), x^k-P_C(z^k)\ra\right\},\;
x^{k,j}=\theta^{j}P_C(z^k)+(1-\theta^{j})x^k.
\]
\item [{\bf(d)} ] \emph{Exogenous stepsize before projecting}: $\alpha_k =
1$ for all $k$ and $\beta_k$ given
by
\begin{equation}\label{(4)}
\beta_k=\frac{\delta_k}{\|\nabla f(x^k)\|}\,,\; \mbox{ with } \;
\sum_{k=0}^\infty\delta_k=\infty \; \mbox{ and }\;
\sum_{k=0}^\infty\delta_k^2<\infty.
\end{equation}
\end{retraitsimple}
Strategy (a) was analyzed in \cite{goldtein} and its weak convergence was proved under Lipschitz continuity of $\nabla f$. The
main difficulty is the necessity of taking $\beta \in (0,2/L)$, where $L$ is the Lipschitz constant of $\nabla f$, which is in general unknown; see also \cite{bertsekas}.

Note that Strategy (b) requires one projection onto $C$ for each
step of the inner loop resulting from the Armijo line search.
Therefore, many projections might be performed for each iteration
$k$, making Strategy (b) inefficient when the projection onto $C$ is
not explicitly computed. On the other hand, Strategy (c) demands
only one projection for each outer step, i.e., for each iteration
$k$. Strategies (b) and (c) are the constrained versions of the
line search proposed in \cite{Armijo-1} for solving unconstrained
optimization problems. Under existence of minimizers and convexity
assumptions for problem~\eqref{1}, it is possible to prove, for
Strategies (b) and (c), convergence of the whole sequence to a
minimizer of $f$ in finite dimensional spaces; see, for instance, \cite{burachik,yun-luis}.
No additional assumption on boundedness of level sets is required,
as shown in \cite{IUSEM-2003}.

Strategy (d), as its counterpart in the unconstrained case, fails to
be a descent method. Furthermore, it is easy to show that this
approach satisfies $\|x^{k+1}-x^{k}\|\le \delta_k$ for all $k$, with
$\delta_k$ exogenous and satisfying \eqref{(4)}. This reveals that
convergence of the sequence of points generated by this approach can
be very slow: in view of \eqref{(4)}, stepsizes are {\em small}
(notice that Strategies (b) and (c) allow for occasionally long
stepsizes because both strategies employ all information available at each iteration).
Moreover, Strategy (d) does not take into account functional values for determining the stepsizes.  These
characteristics, in general, entail poor computational performance.
The strategy's redeeming feature is that its convergence properties
also hold in the nonsmooth case, in which the two Armijo line
searches given by (b) and (c) may be unsuccessful; see
\cite{polyak_libro}. By assuming existence of solutions of
problem~\eqref{1}, defined in an arbitrary Hilbert space $\HH$, and
replacing, at each iteration, $\nabla f(x^k)$ by any subgradient $s_k$ of $f$ at $x^k$,
the work \cite{AIS} establishes that the sequence $(x^k)_{k\in\NN}$
generated by Strategy (d) converges weakly to a solution of
problem~\eqref{1}, providing that the subdifferential of $f$ is
bounded on bounded sets.

In finite dimensional spaces and without assuming convexity of the
function $f$, convergence results for the above strategies closely
mirror the ones for the steepest descent method in the unconstrained
case: cluster points may fail to exist, even when \eqref{1} has
solutions. However, if cluster points exist, they are stationary and
feasible; see for instance \cite[Section 2.3.2]{bertsekas}. The work
\cite{bertsekas-2} proves convergence of the sequence of points
generated by Strategy (b) to a stationary point of problem~\eqref{1}
by assuming  that the starting iterate $x^0$ belongs to a bounded
level set of $f$.

As already mentioned, in this paper we are interested in weak and
strong convergence of projected gradient methods applied to convex
programs as~\eqref{1}. To the best of our knowledge, weak
convergence of the projected gradient method has only been shown
under the assumption of  Lipschitz continuity of $\nabla f$ or using
exogenous stepsize, like Strategy (d) above. In the present work we
prove,  without Lipschitz continuity assumption, weak convergence of
the projected gradient method employing Strategy (c). Moreover, we
propose a few modifications of Strategy (c) in order to ensure that
the resulting method is strongly convergent.

The paper is organized as follows. The next section provides some
definitions and preliminary results that will be used in the
remainder of this work. Weak convergence of the projected gradient
methods is presented in Section $3$. In Section $4$ we propose and
study strong convergence of a modified projected gradient method.
Finally, some comments and remarks are presented in Section $5$.

Our notation is standard: the inner product in $\HH$ is denoted by
$\la\cdot,\cdot\ra$ and the norm induced by the inner product is
$\|\cdot\|$. For an element $x\in \HH$, we define the metric or nearest point
projection of $x$ onto $C$, denoted by $P_C(x)$, as the unique point
in $C$ such that $\|P_C(x)-y\|\le\|x-y\|$ for all $y\in C$. The
indicator function of $C$, written as $I_C$, is given by $I_C(x) =0$
if $x \in C$, and $I_C(x) =\infty$  otherwise, and the normal cone
to $C$ is $\mathcal{N}_C=\partial I_C$. Furthermore, if we define
the function $\hat{f}:=f+I_C$, then problem \eqref{1} is
equivalent to find $x\in \HH$ such that $0\in
\partial \hat{f}(x)=
\partial f(x)+\mathcal{N}_C(x)$, which will be used in the proof of Theorem \ref{todos_ptos_de_acumulacion_estan_S_*}.
\section{Preliminaries}

We begin by stating the (one-sided) directional derivative of $f$ at $x \in \mbox{dom} (f)$ in the direction $d$, that is
\[
f^\prime(x;d) := \lim_{t \to 0^+} \frac{f(x + td) - f(x)}{t},
\]
when the limit exists. From here on we drop the adjective ``one-sided'' and refer to this simply as the directional derivative. If the directional derivative $f'(x,d)$ exists for all directions $d$ and the functional $\nabla f(x): \HH\rightarrow \RR$
defined by $\langle \nabla
f(x),\cdot\,\rangle:=f'(x;\cdot\,)$ is linear and bounded, then we
say that $f$ is Gateaux differentiable at $x$, and $\nabla f(x)$
is called the Gateaux derivative. Every convex and  lower semicontinuous function $f:\HH \to \RR$
that is Gateaux differentiable
at $x$ is also continuous at $x$.

We now recall some necessary and sufficient optimality
conditions for problem~\eqref{1}, whose proof can be found in
\cite[Prop. 17.4]{librobauch}.
\begin{proposition} \label{1.4.3}
Let $f \colon \HH \rightarrow \RR$ be a proper convex and Gateaux
differentiable function. Then the point $x_*\in C$ is a minimizer of
problem \eqref{1} if and only if $0 \in \nabla f(x_*)+\mathcal{N}_C(x_*)$ if
and only if $\la\nabla f(x_*),y-x_*\ra\ge0$ for all $y\in C$ if and
only if $x_*=P_C(x_*-\beta\nabla f(x_*))$ with $\beta>0$.
\end{proposition}
We state some well known facts about the metric projection, first appeared in \cite{nashed}; see also Lemmas    $1.1$    and    $1.2$    in    \cite{zarantonelo}.
\begin{proposition}\label{popiedades_projeccion}
Let $K$  be a nonempty, closed and convex set in $\HH$. For all
$x,y\in \HH$ and all $z\in K$, the following properties hold:
\item [ {\bf(i)} ] $\la x-P_K(x),
z-P_K(x)\ra\leq0$.
\item [ {\bf(ii)} ] $\la z-y,z-P_K(y)\ra \geq \|
z-P_K(y)\|^2$.
\end{proposition}
\noindent Next we deal with the so called quasi-Fej\'er convergence and its
properties.
\begin{definition}\label{def-cuasi-fejer}
Let $S$ be a nonempty subset of $\HH$. A sequence $(x^k)_{k\in\NN}$
in $\HH$ is said to be quasi-Fej\'er convergent to $S$ if and only
if there exists a summable sequence $(\epsilon_k)_{k\in \NN}$ such
that for all $x \in S$, $\| x^{k+1}-x\|^2 \leq \|
x^{k}-x\|^2+\epsilon_k$ for all $k$.
\end{definition}
This definition originates in \cite{browder} and has been further
elaborated in \cite{IST,comb}. A useful result on quasi-Fej\'er sequences is the following, which is proven in \cite[Lemma 6]{browder}.
\begin{lemma}\label{cuasi-Fejer}
If $(x^k)_{k\in\NN}$ is quasi-Fej\'er convergent to $S$, then:
\item [ {\bf(i)} ] $(x^k)_{k\in\NN}$ is bounded;
\item [ {\bf(ii)} ] if all weak cluster
point of $(x^k)_{k\in\NN}$ belong to $S$, then the sequence
$(x^k)_{k\in\NN}$ is weakly convergent.
\end{lemma}
We finalize this section by showing that uniform continuity of the derivative $\nabla f$ is a weaker assumption than Lipschitz continuity of $\nabla f$ in $\HH$.
\begin{example}\label{ex}
Take in problem~\eqref{1} $f(x)=(1/p)\|x\|^p$, with $p>1$. It is not difficult to show that $f$ is a convex function and that $\nabla f$ is
uniformly continuous for all $p>1$. However, $\nabla f$ is globally
Lipschitz continuous only for $p = 2$.
\end{example}

\section{Weak Convergence of the Projected Gradient Method}\label{gradiente-projectado}

In this section we state the classical projected gradient method,
with the line search along the feasible direction (i.e., Strategy (c)).
Provided that the
underlying problem has a solution, we show that the sequence of
points generated by the project gradient method converges weakly to a
solution of the problem. We will not assume Lipschitz continuity of the mapping $\nabla f$.

We  now recall the formal definition of the projected gradient
method. Let $(\beta_k)_{k\in \NN}$ be a sequence such that
$(\beta_k)_{k\in \NN}\subset [\check{\beta},\hat{\beta}] $ with
$0<\check{\beta} \leq \hat{\beta}<\infty$, and be $\theta,
\delta\in(0,1)$. The algorithm is stated as follows:
\begin{center}\fbox{\begin{minipage}[b]{\textwidth}
\noindent {\bf Algorithm A1}

\medskip

\noindent{\bf Initialization Step:} Take $x^0\in C$ and set $k=0$.

\noindent{\bf Iterative Step 1:} Given $x^k$, compute
\begin{equation}\label{paso1_*}
z^k:=x^k-\beta_k \nabla f(x^k)\;.
\end{equation}
\noindent{\bf Stop Criterion:} If $x^k=P_C(z^k)$ then stop.

\noindent {\bf Inner Loop:}  Otherwise, set $\alpha_k=\theta^{j(k)}$, where
\begin{equation}\label{12*}
j(k):=\min \left\{j\in \NN \mid f(x^{k,j})\le
f(x^k)-\delta\theta^{j}\la\nabla f(x^k), x^k-P_C(z^k)\ra\right\},\; x^{k,j}=\theta^{j}P_C(z^k)+(1-\theta^{j})x^k.
\end{equation}
\noindent{\bf Iterative Step 2:}
Compute
\begin{equation}\label{paso3_*}
x^{k+1}= \alpha_k P_C(z^k)+(1-\alpha_k)x^k.
\end{equation}
Set $k=k+1$ and go back to {\bf Iterative Step 1}. 
\end{minipage}}\end{center}
\noindent It follows from Proposition~\ref{popiedades_projeccion}(ii) that the
iterates of {\bf Algorithm A1} satisfy
\begin{equation}\label{innergrad}
\la\nabla f(x^k), x^k-P_C(z^k )\ra\geq \frac{1}{\beta_k}\|x^k-P_C(z^k)\|^2\;\mbox{ for all $k$}.
\end{equation}
Moreover, if {\bf Algorithm A1} stops  then $x^k=P_C(z^k) =
P_C(x^k-\beta_k\nabla f(x^k))$. Since $\beta_k\geq \check \beta>0$,
it follows from Proposition~\ref{1.4.3} that $x^k$ is a solution to
problem~\eqref{1}. Moreover, from~\eqref{12*} and \eqref{innergrad}
we have
\[
f(x^{k+1}) \leq f(x^k)-\delta\alpha_k \la\nabla f(x^k),x^k-P_C(z^k)\ra
\leq f(x^k)- \delta\frac{\alpha_k}{\beta_k}\|x^k-P_C(z^k)\|^2\;\mbox{ for all $k$}.
\]
Therefore, if the algorithm does not stop we obtain the inequality
\begin{equation*}
\delta\frac{\alpha_k}{\hat \beta}\|x^k-P_C(z^k)\|^2\leq f(x^k)-f(x^{k+1}),
\end{equation*}
showing that $(f(x^k))_{k \in \NN}$ is a monotone decreasing sequence. Since such sequence is bounded from below by the optimal value of problem~\eqref{1},
we conclude that $\lim_{k \to \infty }(f(x^k)-f(x^{k+1}))=0$. It thus follows from the above inequality that
\begin{equation}\label{limite}\lim_{k\to\infty} \alpha_k\,\| x^k-P_C(z^k)\|^2 = 0,
\end{equation}
a crucial result to be considered in Theorem~\ref{todos_ptos_de_acumulacion_estan_S_*} below.
 In the following we show that the
\emph{inner loop} in {\bf Algorithm A1} is well-defined.
\begin{proposition}\label{propdef}
Let $k$ be a fixed iteration counter, and suppose that $x^k$ is not
a solution to problem~\eqref{1}. Then, after finitely many steps the
\emph{inner loop} in {\bf Algorithm A1} finds $\alpha_k=\theta^{j(k)}$
satisfying~\eqref{12*}.
\end{proposition}
\begin{proof} The proof is by contradiction: suppose
that~\eqref{12*} does not hold for all $j\ge0$, i.e.,
\[
\frac{f(x^k+\theta^{j}(P_C(z^k)-x^k))- f(x^k)}{\theta^j}> -\delta
\la \nabla f(x^k), x^k-P_C(z^k)\ra \; \mbox{ for all }\; j\geq 0.
\]
Passing to the limit when $j$ goes to infinity and using the
Gateaux differentiability of $f$ and \eqref{innergrad}, we
conclude that
\[
0 \geq(1-\delta) \la \nabla f(x^k), x^k-P_C(z^k)\ra\geq
\frac{(1-\delta)}{\beta_k}\|x^k-P_C(z^k)\|^2\geq
\frac{(1-\delta)}{\hat{\beta}}\|x^k-P_C(z^k)\|^2,
\]
which contradicts, by Proposition~\ref{1.4.3}, the assumption that $x^k$ is not a solution to the problem.
\end{proof}
\noindent To the best of our knowledge, from now on all the presented results
are new in Hilbert spaces.
\begin{lemma}\label{qf}
The sequence generated by {\bf Algorithm A1} is quasi-Fej\'er convergent
to $S_*$.
\end{lemma}
\begin{proof}
Take any $x_*\in S_*$. Note that $\|x^{k+1}-x^k\|^2+\|x^{k}-x_*\|^2
-\|x^{k+1}-x_*\|^2=2\la x^{k}-x^{k+1},x^k-x_*\ra $. Moreover,
\begin{align*} 2\la x^{k}-x^{k+1},x^k-x_*\ra=&\, 2\alpha_k\la
P_C(z^{k})-x^{k},x_*-x^k\ra\\[1ex]
=&\, 2\alpha_k\beta_k\la \nabla f(x^k),x^k-x_*\ra -2\alpha_k\la P_C(z^{k})-x^{k}+\beta_k\nabla f(x^k),x^k-x_*\ra \\[1ex]
=&\,2\alpha_k\beta_k\la \nabla f(x^k),x^k-x_*\ra -2\alpha_k\la P_C(z^{k})-x^{k}+\beta_k\nabla f(x^k),x^k-P_C(z^k)\ra \\
&\,-2\alpha_k\la P_C(z^{k})-(x^{k}-\beta_k\nabla f(x^k)),P_C(z^k)-x^*\ra \\[1ex]
\geq&\,2\alpha_k\beta_k\la \nabla f(x^k),x^k-x_*\ra -2\alpha_k\la P_C(z^{k})-x^{k}+\beta_k\nabla f(x^k),x^k-P_C(z^k)\ra \\[1ex]
\geq&\,2\alpha_k\beta_k(f(x^k)-f(x_*)) -2\alpha_k\la P_C(z^{k})-x^{k}+\beta_k\nabla f(x^k),x^k-P_C(z^k)\ra \\[1ex]
\geq&\,-2\alpha_k\la P_C(z^{k})-x^{k}+\beta_k\nabla f(x^k),x^k-P_C(z^k)\ra \\[1ex]
  = &\, 2\alpha_k\|x^k-P_C(z^k)\|^2 - 2\alpha_k\beta_k \la \nabla f(x^k),x^k-P_C(z^k)\ra
,
\end{align*}
where the first equality follows from~\eqref{paso3_*}, the second
one by adding and subtracting $\beta_k\nabla f(x^k)$, and the third
equality follows from adding and subtracting $P_C(z^k)$. The first
inequality above is due to
Proposition~\ref{popiedades_projeccion}(i), the second one follows
from convexity, and the third inequality holds because $x_*$ is a
solution to problem~\eqref{1}. We thus have shown that
\[\|x^{k+1}-x_*\|^2\le\|x^{k}-x_*\|^2+\|x^{k+1}-x^k\|^2 -2\alpha_k \|x^k-P_C(z^k)\|^2
+2\alpha_k\beta_k \la \nabla
f(x^k),x^k-P_C(z^k)\ra.\]  Since $x^{k+1}-x^k=\alpha_k(P_C(z^k) -x^k)$ by \eqref{paso3_*}, we conclude that
\[
\|x^{k+1}-x_*\|^2\le\|x^{k}-x_*\|^2+\alpha_k^2\|x^k-P_C(z^k)\|^2 -2\alpha_k \|x^k-P_C(z^k)\|^2
+2\alpha_k\beta_k \la \nabla
f(x^k),x^k-P_C(z^k)\ra.
\]
Moreover, we have that $\alpha_k^2 - 2\alpha_k\leq -\alpha_k$,
because $0\leq \alpha_k\leq 1$ in the \emph{inner loop} of {\bf Algorithm
A1}. This gives
\[
\|x^{k+1}-x_*\|^2\le\|x^{k}-x_*\|^2-\alpha_k \|x^k-P_C(z^k)\|^2
+2\frac{\hat{\beta}}{\delta}\left(f(x^k)-f(x^{k+1})\right),
\]
where we have used \eqref{12*}.
In order to show that $(x^k)_{k\in\NN}$ is
quasi-Fej\'er convergent to $S_*$, it remains to prove that $$\epsilon_k:=-\alpha_k\|x^k-P_C(z^k)\|^2+2\frac{\hat{\beta}}{\delta}\left(f(x^k)-f(x^{k+1})\right)$$
forms a convergent series. Indeed, this is true by the following development
\[
\sum_{k=0}^\infty\epsilon_k\le 2\frac{\hat{\beta}}{\delta}\sum_{k=0}^\infty\left(f(x^k)-f(x^{k+1})\right)
= 2\frac{\hat{\beta}}{\delta}\left(f(x^0)-\lim_{k\to\infty} f(x^{k+1})\right)
\leq 2\frac{\hat{\beta}}{\delta}\left(f(x^0)-f(x_*)\right)< \infty.
\]\end{proof}
\noindent We are now ready to prove the main result of this section.
\begin{theorem}\label{todos_ptos_de_acumulacion_estan_S_*}
Suppose that $\nabla f$ is uniformly continuous on bounded sets.
Then the sequence $(x^k)_{k\in\NN}$ generated by {\bf Algorithm A1} is
bounded and each of its weak cluster points belongs to $S_*$.
\end{theorem}
\begin{proof}
Since the sequence $(x^k)_{k\in\NN}$ is quasi-Fej\'er convergent to $S_*$, it is
bounded. Therefore, there exists a subsequence
$(x^{i_k})_{k\in \NN}$ of $(x^k)_{k\in\NN}$ that converges weakly to
some $x_*$.
Moreover, since $\nabla f$ is uniformly continuous on bounded sets we conclude that $(\nabla f(x^k))_{k \in \NN}$ is also a bounded sequence.
Thus, it follows from~\eqref{paso1_*} and \eqref{12*} that $(P_C(z^k))_{k\in \NN}$ is a bounded sequence as well.

\noindent We now split our analysis into two distinct cases.

\medskip

\noindent {\bf Case 1.} Suppose that the sequence $(\alpha_k)_{k\in \NN}$ does
not converge to $0$, i.e. there exists a subsequence
$(\alpha_{i_k})_{k\in \NN}$ of $(\alpha_{k})_{k\in \NN}$ and some
$\alpha>0$ such that $\alpha_{i_k}\geq\alpha$ for all $k$. Let $w^k:=P_C(z^{k})$; it follows from \eqref{limite} and our assumption on $\nabla f$ that
\begin{equation}\label{limite1_DI}
\dsty\lim_{k\rightarrow\infty}\,\| x^{i_k}-w^{i_k}\| = 0\;\mbox{ and }\; \dsty\lim_{k\rightarrow\infty}\| \nabla f(x^{i_k})-\nabla f(w^{i_k})\| =0.
\end{equation}
Let $x_*$ be a weak cluster point of the subsequence
$(x^{i_k})_{k\in \NN}$. By \eqref{limite1_DI}, it is also a weak
cluster point of $(w^{i_k})_{k\in \NN}$. Without loss of generality,
we assume that $(x^{i_k})_{k\in \NN}$ and $(w^{i_k})_{k\in \NN}$
converges weakly to $x_*$. In order to prove that $x_* \in S_*$, we define the function $\hat{f}:=f+I_C$. It is well known that
$
\partial \hat{f}(x):= \nabla f(x)+\mathcal{N}_C(x)
$, for all $x\in C$,  is a maximal monotone, and that $0\in
\partial\hat{f}(x)$ if and only if $x\in S_*$; see for instance
\cite{Rock}. Therefore, we need to show that $0 \in \partial \hat
f(x_*)$. In order to do that, we take $(x,u)\in G(\partial \hat{f})$
with $x\in C$. Thus, $u\in\partial \hat{f}(x)=\nabla
f(x)+\mathcal{N}_C(x)$, implying that $u-\nabla f(x)\in
\mathcal{N}_C(x)$. So, we have
$
\la x-y,u-\nabla f(x) \ra\geq0$ for all $y\in C
$. 
In particular,
\begin{equation}\label{*}
\la x-w^{i_k},u\ra\geq \la x -w^{i_k}, \nabla f(x)\ra.
\end{equation}
On the other hand, since $ w^{i_k}=P_C(x^{i_k}-\beta_{i_k} \nabla
f(x^{i_k}))$ and $x^{i_k}\in C$, it follows from Proposition
\ref{popiedades_projeccion}(i), with $K=C$ and
$x=x^{i_k}-\beta_{i_k}\nabla f(x^{i_k})$, that $\la
x-w^{i_k},x^{i_k}-\beta_{i_k} \nabla f(x^{i_k})-w^{i_k} \ra\leq0$
for all $x\in C$ and ${i_k}\geq0$. Rearranging terms and taking into
account that $\beta_{i_k}>0$, we get
\[
\left\la x-w^{i_k},\frac{x^{i_k}-w^{i_k}}{\beta_{i_k}}-\nabla f(x^{i_k}) \right\ra\leq0
\quad\forall\, x\in C\,\,\mbox{and}\,\,{i_k}\geq0.
\]
Together with \eqref{*}, we conclude that
\begin{align}\nonumber
\dsty\la x-w^{i_k},u\ra&\geq\la x-w^{i_k}, \nabla f(x) \ra\geq\la x-w^{i_k}, \nabla f(x)
\ra+\left\la x-w^{i_k}, \frac{x^{i_k}-w^{i_k}}{\beta_{i_k}} - \nabla f(x^{i_k})
\right\ra\\\nonumber
&=\la x-w^{i_k}, \nabla f(x)-\nabla f(w^{i_k}) \ra + \la x-w^{i_k},
\nabla f(w^{i_k})-\nabla f(x^{i_k})\ra+\left\la x-w^{i_k},
\frac{x^{i_k}-w^{i_k}}{\beta_{i_k}} \right\ra.
\nonumber
\end{align}
Monotonicity of $\nabla f$ gives $\la x-w^{i_k}, \nabla f(x)-\nabla
f(w^{i_k}) \ra \geq 0$. Thus,
\begin{align}\nonumber
\dsty\la x-w^{i_k},u\ra&\geq\la x-w^{i_k}, \nabla f(w^{i_k})-\nabla f(x^{i_k}) \ra+\left\la x-w^{i_k},
\frac{x^{i_k}-w^{i_k}}{\beta_{i_k}} \right\ra\\\nonumber
&\geq-\| x-w^{i_k}\| \left(\| \nabla f(w^{i_k})-\nabla f(x^{i_k})\| +\frac{1}{\beta_{i_k}}\| w^{i_k}-x^{i_k}\| \right)\\
&\geq -\| x-w^{i_k}\| \left(\| \nabla f(w^{i_k})-\nabla f(x^{i_k})\|
+\frac{1}{\check{\beta}}\| w^{i_k}-x^{i_k}\| \right),\nonumber
\end{align}
where we have used Cauchy-Schwarz inequality in the third
inequality and the fact that $\beta_k\geq\check{\beta}>0$ for all
$k$ in the last one. Remember that $\{w^{i_k}\}_{k \in\NN}$ is
bounded and converges weakly to $x_*$. Thus, passing to the limit in
the above relations and using~\eqref{limite1_DI}, we obtain
\[
\la x-x_*,u \ra\geq0 \quad\quad\forall\,(x,u)\in G(\partial
\hat{f}).
\]
Since $\partial \hat{f}$ is maximal monotone, it follows from the above inequality that $(x_*,0)\in G(\partial \hat{f})$, implying
that $0\in
\partial \hat{f}(x_*)=\nabla f(x_*)+\mathcal{N}_C(x_*)$ and hence $x_*\in S_*$.

\medskip

\noindent{\bf Case 2.} Suppose now that
$\lim_{k\rightarrow\infty}\alpha_{k}=0$. Take, with $\alpha_k>0$,
\begin{equation}\label{2*paso}
\hat{y}^k=\frac{\alpha_k}{\theta}P_C(z^k)+\left(1-\frac{\alpha_k}{\theta}\right)x^k=x^k-\frac{\alpha_k}{\theta} (x^k-w^k).
\end{equation}
It follows from the definition of $j(k)$ in \eqref{12*} that
$
f(\hat{y}^k)-f(x^k)>-\delta \frac{\alpha_k}{\theta}\la\nabla f(x^k), x^k-w^k\ra
$. 
Thus,
\begin{align*}\label{eq-garantiza-xk-Pcxk-va-cero}
\delta\frac{\alpha_k}{\theta}\la\nabla f(x^k), x^k-w^k\ra&>f(x^k)-f(\hat{y}^k)\geq\la\nabla f(\hat{y}^k), x^k-\hat{y}^k\ra
=\frac{\alpha_k}{\theta}\la\nabla f(\hat{y}^k), x^k-w^k\ra\\
&=\frac{\alpha_k}{\theta}\dsty\la \nabla f(\dsty \hat{y}^k)-\nabla f(x^k),x^{k}-w^k\ra+\frac{\alpha_k}{\theta}\dsty\la
\nabla f(x^k),x^{k}-w^k\ra\\
&\geq -\frac{\alpha_k}{\theta}\| \nabla f(\hat{y}^k)-\nabla f(x^k)\| \| x^k-w^k\|
+\frac{\alpha_k}{\theta}\dsty\la
\nabla f(x^k),x^{k}-w^k\ra,
\end{align*}
where we have used convexity of $f$ in the second inequality and
Cauchy-Schwarz inequality in the last one. Rearrangement of terms and
using \eqref{innergrad} yield
$$
\| \nabla f(\hat{y}^k)-\nabla f(x^k)\| \| x^k-w^k\|
\geq(1-\delta)\la\nabla f(x^k), x^k-w^k\ra\ge
\frac{(1-\delta)}{\beta_k}\|x^k-w^k\|^2,
$$ which implies
\begin{equation}\label{aux0}
\| \nabla f(\hat{y}^k)-\nabla f(x^k)\|\ge
\frac{(1-\delta)}{\hat{\beta}}\| x^k-w^k\|.
\end{equation}

Since both sequences $(x^k)_{k\in \NN}$ and $(w^k)_{k \in \NN}$ are
bounded and $\lim_{k\to \infty} \alpha_k =0$, it follows from
\eqref{2*paso} that $\lim_{k\to \infty} \|\hat y^k -x^k\|=0$. As
$\nabla f$ is uniformly continuous on bounded sets, we get
$\lim_{k\rightarrow\infty}\| \nabla f(\hat{y}^k)-\nabla
f(x^k)\| =0$. It thus follows from \eqref{aux0} that
$\lim_{k\rightarrow\infty}\,\|
x^k-w^k\| =0.$
We have show that {\bf Case 2} also satisfies the key relations in \eqref{limite1_DI} of {\bf Case 1}.
Hence, the remain of the proof can be done similarly to {\bf Case 1}, \emph{mutatis mutandis}.
\end{proof}

\begin{theorem}
Suppose that $\nabla f$ is uniformly continuous on bounded sets.
Then the sequence $(x^k)_{k\in\NN}$ generated by {\bf Algorithm A1}
converges weakly to a solution of problem~\eqref{1}.
\end{theorem}
\begin{proof}
By Lemma \ref{qf} $(x^k)_{k\in\NN}$ is quasi-Fej\'er convergent to
$S_*$ and by Theorem \ref{todos_ptos_de_acumulacion_estan_S_*} all
weak cluster points of $(x^k)_{k\in\NN}$ belong to $S_*$. The result thus follows from  Lemma~\ref{cuasi-Fejer}(ii).
\end{proof}

\section{Strongly Convergent Projected Gradient Method}

In this section we consider a modification of the projected gradient method forcing strong convergence in Hilbert spaces.
The modified projected method, employing line search (c), was inspired by Polyak's method \cite{poljak, yunier-iusem-4, yunier-wlo} for nondifferentiable optimization. The method uses a similar idea to that exposed in \cite{yunier-iusem-1, SOLODOV_SVAITER}, with the same goal, upgrading weak to strong convergence.

Additionally, our algorithm has the distinctive feature that the
limit of the generated sequence is the solution of the
problem closest to the initial iterate $x^0$. This property is useful in
many specific applications, e.g. in image reconstruction \cite{ImageR-1,ImageR-2,ImageR-3} and in minimal norm solution problems,
as discussed in \cite{wim-yunier-wlo}. We
emphasize that this feature is of interest also in finite-dimensional spaces, differently from the strong versus weak convergence
issue.

\subsection{Some Comments on Strong Convergence}

Clearly weak and strong convergence are only distinguishable in the
infinite-dimensional setting. Naturally, even when we have
to solve infinite-dimensional problems, numerical implementations of
algorithms are performed in finite-dimensional
approximations of these problems. Nevertheless, it is interesting to
have good convergence theory for the infinite-dimensional setting in order
to guarantee robustness and stability of the finite-dimensional approximations. This issue is closely related to the so-called Mesh
\emph{Independence Principle}  presented in \cite{[2],[1],[15]}. This principle relies on infinite-dimensional convergence to predict the convergence
properties of a discretized finite-dimensional method. Moreover, the
Mesh Independence Principle provides theoretical justification for
the design of refinement strategies, which are crucial for having appropriate approximation to the true solution of the infinite-dimensional problem being solved.
We suggest the reader to
see \cite{[13]}, where a variety of applications are described. A
strong convergence principle in Hilbert spaces is extensively
analyzed in \cite{bauschke, alves-jeff}.

The importance of strong convergence is also underlined in
\cite{guler}, where it is shown, for the proximal-point algorithm,  that the rate of convergence of the value sequence $(f(x^k))_{k\in \NN}$ is better when
$(x^k)_{k\in\NN}$ converges strongly. It is important to say that only weak
convergence has been established for the projected gradient method
in Hilbert spaces; see \cite{5, goldtein}. In these cases the weak
convergence has been established by assuming Lipschitz continuity of
$\nabla f$ or by employing exogenous stepsizes that may lead to small-length steps, as discussed in the Introduction. In our scheme we use
the classical Armijo line search along the feasible direction
establishing the strong convergence.

\subsection{Algorithm and Convergence Analysis}

Let $(\beta_k)_{k\in \NN}$ be a sequence such that $(\beta_k)_{k\in
\NN}\subset [\check{\beta},\hat{\beta}] $ with $0<\check{\beta} \leq
\hat{\beta}<\infty$, and be $\theta, \delta\in(0,1)$. The algorithm
of the proposed strongly convergence projected gradient is stated as
follows.
\begin{center}\fbox{\begin{minipage}[b]{\textwidth}
\noindent {\bf Algorithm A2}

\medskip

\noindent{\bf Initialization Step:} Take $x^0\in C$ and set $f^{lev}_{-1} =\infty$.

\noindent{\bf Iterative Step 1:} Given $x^k$, compute $
z^k=x^k-\beta_k \nabla f(x^k).$

\noindent{\bf Stop Criterion 1:} If $x^k=P_C(z^k)$ then stop.

\noindent {\bf Inner Loop:}  Find $j(k)$ as in \eqref{12*},
set $\alpha_k=\theta^{j(k)}$ and $f^{lev}_k =\min\{f^{lev}_{k-1},f( x^{k,j(k)})\}$.

\noindent{\bf Iterative Step 2:}
Define
\begin{equation*}\label{Hk}
H_k:=\left\{x\in \HH\mid \la \nabla f(x^k),x-x^k\ra+f(x^k)-f^{lev}_k\leq 0\right\},\; W_k:=\left\{x\in \HH\mid \la x-x^k,x^0-x^k\ra\leq 0\right\}.
\end{equation*}
Compute
\begin{equation}\label{paso3_**}
\dsty x^{k+1}:=P_{C\cap W_k\cap H_k}(x^0).
\end{equation}
\noindent{\bf Stop Criterion 2:} If $x^{k+1}=x^k$ then stop.

\noindent Otherwise, set $k=k+1$ and go back to {\bf Iterative Step 1}.
\end{minipage}}\end{center}
\noindent {\bf Algorithm A2} is a particular case of Algorithm 2 in
\cite{yunier-wlo}, for nonsmooth convex optimization. In contrast to
the algorithm in \cite{yunier-wlo}, {\bf Algorithm A2} gives a
straightforward way to define the level sequence
$(f^{lev}_k)_{k \in\NN}$.

Suppose that $x^k \notin S_*$. Since
$\la\nabla f(x^k), x^k-P_C(z^k)\ra \geq
\frac{1}{\beta_k}\|x^k-P_C(z^k)\|$, the definition of $f^{lev}_k$
does satisfies the inequalities given in \cite[Eq. (4)]{yunier-wlo}:
\begin{equation}\label{>flev}
f(x^k)> f^{lev}_k \geq f_* \; \mbox{ for all }\; k\geq 0.
\end{equation}
We thus conclude that if $x^k \notin S_*$ then $f(x^k)> f^{lev}_k$,
yielding that $x^k \notin H_k$. In order to analyze convergence
of {\bf Algorithm A2} we present below some key inequalities.
\begin{lemma}
For all $k\ge 0$  it holds that
\begin{equation}\label{des-1}
\| x^{k+1}-x^0\| ^2\geq\| x^k-x^0\|^2+ \|x^k-x^{k+1}\|^2,\;\; \mbox{with }\;
\end{equation}
\begin{equation}\label{des-2} \| x^k-x^{k+1}\| \geq \frac{f(x^k)-f^{lev}_k}{\|\nabla f(x^k)\|}\geq\delta\frac{\alpha_k}{\hat{\beta}}\,
\frac{\| x^k-P_C(z^k)\|^2}{\|\nabla f(x^k)\|}\geq 0.
\end{equation}
\end{lemma}
\begin{proof}
\noindent Since $x^{k+1}\in W_k$, then
$$
0\geq\la
x^{k+1}-x^k,x^0-x^k\ra=\frac{1}{2}\left(\|x^{k+1}-x^k\|^2-\|x^{k+1}-x^0\|^2+\|x^k-x^0\|^2\right),
$$ 
which implies \eqref{des-1}.
Moreover, $x^{k+1}$ belongs to $H_k$ and thus
\[
\la \nabla f(x^k), x^k-x^{k+1}\ra \geq f(x^k)-f^{lev}_k,\;
\mbox{ i.e, }\; \| x^k-x^{k+1}\| \geq \frac{f(x^k)-f^{lev}_k}{\|\nabla f(x^k)\|}.
\]
Using \eqref{12*} and \eqref{innergrad}, $$f(x^k)-f^{lev}_k\geq
f(x^k)-f(x^{k,j(k)})\geq \delta \alpha_k \la\nabla f(x^k),
x^k-P_C(z^k)\ra \geq \delta\frac{\alpha_k}{\beta_k}\,\|
x^k-P_C(z^k)\|^2\geq 0, $$ and the result follows because $\beta_k
\leq \hat{\beta}$.
\end{proof}
\noindent We now put together some results from \cite{yunier-wlo}, yielded by \eqref{paso3_**} and \eqref{>flev}.
\begin{lemma}\label{lemma:wlo}
Assume that $\nabla f$ is uniformly bounded on bounded sets, and
that the solution set $S^*$ of problem \eqref{1} is nonempty. Let
$(x^k)_{k\in\NN}$ be the sequence of points generated by {\bf Algorithm
A2}, and let $x^*$ be the projection of $x^0$ onto $S^*$, i.e.,
$x^*=P_{S_*}(x^0)$. Then,
\item [ {\bf(i)}] $S_*\subseteq H_k \cap W_k\cap C$ for all $k$;
\item [ {\bf(ii)}] $(x^k)_{k\in\NN}$ is contained in the closed ball centered in $( x^0  + x^*)/2$ and with radius $\|x^* - x^0\|/2$;
\item [ {\bf(iii)}] all weak cluster points of $(x^k)_{k\in\NN}$ belongs to $S^*$.
\end{lemma}
\begin{proof}
It follows from our assumptions on problem~\eqref{1} that $x^*$ is
well-defined. Since \eqref{>flev} holds, Equation (11) in
\cite{yunier-wlo} holds for the choice $i=0$ and $j=\infty$ in the
definition of $I_i^j$ therein. Therefore, both items (i) and
(ii) follows from taking $\hat x =x^0$ in \cite[Prop.
2.2]{yunier-wlo}.

As a result of (ii), $(x^k)_{k\in\NN}$ is bounded and has at least
a weak cluster point. It follows from \eqref{des-1} that the
sequence $(\|x^k-x^0\|)_{k\in \NN}$ is nondecreasing and bounded,
hence convergent. Again, by \eqref{des-1}:
$$0\leq\|x^{k+1}-x^k\|^2\leq\|x^{k+1}-x^0\|^2-\|x^k-x^0\|^2,$$ and we
conclude that
\begin{equation}\label{xk+1-xk-va-cero}
\lim_{k\rightarrow\infty}\| x^{k+1}-x^k\|^2 =0.
\end{equation}
Boundedness of $(\nabla f (x^k))_{k\in \NN}$ follows from the
boundedness of $(x^k)_{k\in\NN}$. Thus, \eqref{des-2} and \eqref{xk+1-xk-va-cero} yields
\[ \lim_{k\to\infty} \alpha_k\,\| x^k-P_C(z^k)\|^2 = 0.\,\]
The proof of item (iii) follows now from repeating the proof of
Theorem \ref{todos_ptos_de_acumulacion_estan_S_*}, cases 1 and 2.
\end{proof}

Since items (i)-(iii) above hold, convergence of {\bf Algorithm A2}
follows from \cite[Thm. 3.4]{yunier-wlo}. We repeat the result here
for completeness.

\begin{theorem} Assume that $\nabla f$ is uniformly bounded
on bounded sets. Then, the sequence $(x^k)_{k\in\NN}$ generated by
{\bf Algorithm A2} converges strongly to $x_* =P_{S_*}(x^0)$.
\end{theorem}
\begin{proof}
It follows from the definition of $x^{k+1}$
that
$
\| x^{k+1}-x^0\| \leq\| x-x^0\|$ for all $ x \in
H_k\cap W_k\cap C$.
In particular, $x_*\in
H_k\cap W_k\cap C$ by Lemma~\ref{lemma:wlo}(i). Thus,
\begin{equation}\label{eq1}
\| x^{k}-x^0\| \leq\| x_*-x^0\| \; \mbox{ for all $k$.}
\end{equation}
By items (ii) and (iii) of Lemma~\ref{lemma:wlo}, $(x^k)_{k\in\NN}$ is
bounded and each of its weak cluster points belongs to $S_*$. Let
$\{x^{i_k}\}$ be any weakly convergent subsequence of
$(x^k)_{k\in\NN}$, and let $\hat{x}\in S_*$ be its weak limit.
Observe that
\begin{align*}
\| x^{i_k}-x_*\| ^2 =&\, \| x^{i_k}-x^0-(x_*-x^0)\| ^2\\
=&\, \| x^{i_k}-x^0\| ^2 + \| x_*-x^0\| ^2-2\la x^{i_k}-x^0, x_*-x^0\ra\\
\leq&\,2\| x_*-x^0\| ^2-2\la x^{i_k}-x^0, x_*-x^0\ra,
\end{align*}
where the inequality follows from (\ref{eq1}). By the weak
convergence of $\{x^{i_k}\}$ to $\hat{x}$, we obtain
\begin{equation}\label{eq2}
\limsup_{k\rightarrow\infty}\| x^{i_k}-x_*\| ^2\leq 2(\| x_*-x^0\|
^2-\la\hat{x}- x^0,x_*-x^0\ra).
\end{equation}
Applying Proposition~\ref{popiedades_projeccion}(i), with $K=S_*$, $x=x^0$
and $z=\hat{x}\in S_*$, and taking into account that $x_*$ is the
projection of $x^0$ onto $S_*$, we have that $\la
x^0-x_*,\hat{x}-x_*\ra\leq 0$. This inequality yields
\begin{align*}
0\geq&\,-\la \hat{x}-x_* , x_*-x^0 \ra=-\la \hat{x}-x^0 ,x_*-x^0 \ra - \la x^0-x_* ,x_*-x^0 \ra\\
\geq&\,-\la \hat{x}-x^0 , x_*-x^0 \ra + \| x_*-x^0\| ^2.
\end{align*}
It follows that
$ \la\hat{x}-x^0,x_*-x^0\ra\geq\| x_*-x^0\| ^2$. By combining this last inequality with \eqref{eq2} we conclude that
$(x^{i_k})_{k\in \NN}$ converges strongly to $x_*$. Thus, we have
shown that every weakly convergent subsequence of $(x^k)_{k\in\NN}$
converges strongly to $x_*$. Hence, the whole sequence
$(x^k)_{k\in\NN}$ converges strongly to $x_*\in S_*$.
\end{proof}

\section{Final Remarks}

It is well-known in Hilbert spaces that global
Lipschitz continuity of the derivative $\nabla f$ is sufficient for
providing convergence of the sequence generated by the projected gradient method, since stepsizes are sufficiently small with respect to the Lipschitz constant.
Naturally, small steps may lead to slow convergence, not mentioning that
having gradients globally Lipschitz is a very restricted assumption.

In this work we dealt with weak and strong convergence of  projected gradient methods for convex (Gateaux) differentiable optimization problems. We focused on the classical Armijo line search along the
feasible direction, eliminating thus the undesired small stepsizes.  Moreover, we relaxed the Lipschitz assumption by supposing only uniform continuity of the derivatives, a much weaker assumption as illustrated in Example~\ref{ex}.
Furthermore, we proposed a strongly
convergent variant of the projected gradient method, which has advantages over the classical projected gradient method.

We hope that this study will serve as basis for future research on other more
efficient variants, as well as including sophisticated line searches
on the gradient methods in Hilbert spaces.

\paragraph {Acknowledgments.}
{\small
\noindent The first author was partially supported by CNPq grants 303492/2013-9,
474160/2013-0 and 202677/2013-3 and by project CAPES-MES-CUBA
226/2012. 

\noindent Research done during a postdoctoral visit of the second author at IMPA - Instituto Nacional de Matem\'atica Pura e Aplicada.

\bibliographystyle{plain}

}
\end{document}